\documentclass[12pt]{article}

\usepackage{subcaption}

\usepackage{float}

\usepackage{multirow} 

\usepackage{tikz}

\usetikzlibrary{calc,backgrounds,arrows,matrix}

\usepackage{enumerate}

\usepackage[colorlinks=true, urlcolor=blue,linkcolor=blue, citecolor=blue]{hyperref}

\usepackage{amssymb,amsmath} 

\usepackage{graphicx} 

\usepackage{graphics} 
\usepackage{color} 

\DeclareMathAlphabet{\eufrak}{U}{}{}{} 
\SetMathAlphabet\eufrak{normal}{U}{euf}{m}{n}
\SetMathAlphabet\eufrak{bold}{U}{euf}{b}{n}

\allowdisplaybreaks

 \def\qu{{\mathord{\mathbb Z}}}

 \def\inte{{\mathord{\mathbb R}}}

 \def\inte{{\mathord{\mathbb N}}}

 \def\sZZ{{\rm Z\kern-.45em{}Z}}

 \def\sQQ{{\kern 0.27em \vrule height1.45ex width0.03em depth0em
           \kern-0.30em \rm Q}}
 \def\qu{{\mathchoice
         {\sQQ}
         {\sQQ}
   {\kern 0.225em \vrule height1.05ex width0.025em depth0em \kern-0.25em \rm Q}
   {\kern 0.180em \vrule height0.78ex width0.020em depth0em \kern-0.20em \rm Q}
         }}
 \def\sGG{{\kern 0.27em \vrule height1.45ex width0.03em depth0em
           \kern-0.30em \rm G}}
 \def\gg{{\mathchoice
         {\sGG}
         {\sGG}
   {\kern 0.225em \vrule height1.05ex width0.025em depth0em \kern-0.25em \rm G}
   {\kern 0.180em \vrule height0.78ex width0.020em depth0em \kern-0.20em \rm G}
         }}

 \newtheorem{prop}{Proposition}[section]
 \newtheorem{lemma}[prop]{Lemma}
 
 \newtheorem{corollary}[prop]{Corollary}
 \newtheorem{theorem}[prop]{Theorem}
 \newtheorem{remark}[prop]{Remark}

\numberwithin{equation}{section}

\def\E{\mathop{\hbox{\rm I\kern-0.20em E}}\nolimits}

 \newcounter{hyp}
\usepackage{aeguill}

\newenvironment{Proof}{\removelastskip\par\medskip \noindent{\em Proof.} \rm}{\penalty-20\null\hfill$\square$\par\medbreak}

\def\bprf{\begin{Proof}}
\def\nprf{\end{Proof}}
\def\bdes{\begin{description}}
\def\ndes{\end{description}}

\newtheorem{thm}{Theorem}[section]

\def\bdef{\begin{defn}}
\def\ndef{\end{defn}}
\def\bthm{\begin{thm}}
\def\nthm{\end{thm}}
\def\bprop{\begin{prop}}
\def\nprop{\end{prop}}
\def\brmk{\begin{remark}}
\def\nrmk{\end{remark}}
\def\bexa{\begin{exa}}
\def\nexa{\end{exa}}
\def\blem{\begin{lem}}
\def\nlem{\end{lem}}
\def\bcor{\begin{cor}}
\def\ncor{\end{cor}}
\def\bexe{\begin{exe}}
\def\nexe{\end{exe}}

\def\rit{\Bbb{R}}

\newcommand{\ee}{\mathbb{E}}

\newcommand{\real}{\mathbb{R}}

\newcommand{\retirer}[1]{$ $\newline 
}

\def\E{\mathop{\hbox{\rm I\kern-0.20em E}}\nolimits}

\def\og{\leavevmode\raise.3ex
     \hbox{$\scriptscriptstyle\langle\!\langle$~}}
\def\fg{\leavevmode\raise.3ex
     \hbox{~$\!\scriptscriptstyle\,\rangle\!\rangle$}~}

\title{\Huge
  Integrability and regularity of the flow of 
  stochastic differential equations with jumps} 

\author{
\large 
Jean-Christophe Breton\footnote{
\href{mailto:jean-christophe.breton@univ-rennes1.fr}{jean-christophe.breton@univ-rennes1.fr}
} 
\\ 
\small 
Univ Rennes\\ 
\small CNRS, IRMAR - UMR 6625
\\ 
\small
263 Avenue du G\'en\'eral Leclerc
\\ 
\small F-35000 Rennes, France
\and
 Nicolas Privault\footnote{
\href{mailto:nprivault@ntu.edu.sg}{nprivault@ntu.edu.sg}
}
\\
\small
Division of Mathematical Sciences 
\\ 
\small 
School of Physical and Mathematical Sciences 
\\ 
\small
Nanyang Technological University 
\\ 
\small 
21 Nanyang Link, Singapore 637371 
}

\allowdisplaybreaks

\begin{document}
\maketitle

\baselineskip0.6cm
 
\vspace{-0.6cm}

\begin{abstract}
  We derive sufficient conditions for the differentiability
  of all orders for the flow of stochastic differential
  equations with jumps, and prove related $L^p$-integrability
  results for all orders. 
  Our results extend similar results obtained in
  \cite{kunitalevy} for first order differentiability
  and rely on the Burkholder-Davis-Gundy (BDG) inequality for time-inhomogeneous Poisson random measures on
  $\real_+\times \real$, for which we provide a new proof.
\end{abstract}

  

\baselineskip0.7cm

\section{Introduction}
\noindent 
In this paper we consider the regularity and integrability
of all orders of the flow of
Stochastic Differential Equations (SDEs) with jumps of the form 
\begin{equation}
\label{eq:edsX}
dX_t(x)=r(t,X_t(x))dt+\sigma(t,X_t(x)) dW_t+\int_{-\infty}^\infty g(t,X_{t^-}(x),y) \big(N(dt,dy)-\nu_t(dy)dt \big), 
\end{equation}
with $X_0(x)=x$, 
\textcolor{black}{where $X_{t^-}$ stands for $\lim_{s\nearrow t} X_s$,}
$g:\real_+ \times\real\times\real\to\real$ is a deterministic measurable function and 
$(W_t)_{t\in \real_+}$, $N(dt,dy)$ are
a standard Brownian motion 
and a Poisson random measure on $\real^+\times\real$ with compensator $\nu_t(dy)dt$, generating a filtration $({\cal F}_t)_{t\geq 0}$. 
\\

\noindent
In the diffusion case,  
the smoothness of the solution flow $x\mapsto X_t(x)$ of stochastic differential equations
 of the form  
\begin{equation}
\label{eq:edsX0}
dX_t(x)=r\big(t,X_t(x)\big)\ dt+\sigma\big(t,X_t(x)\big)\ dW_t
\end{equation} 
with $X_0(x)=x$, where $r:\real_+ \times\real\to\real$, $\sigma: \real_+ \times\real\to\real$
are deterministic coefficients, 
has been studied in \cite{kunitasaintflour},
\cite{kunita}.
In Theorem~II.4.4 of \cite{kunitasaintflour} and Theorem~4.6.5 of \cite{kunita}
it is shown that 
$x\mapsto X_t(x)$ is $k$ times continuously differentiable
when the SDE coefficients of \eqref{eq:edsX0}
are ${\cal C}^k$ functions with globally Lipschitz derivatives.
Such results have been proved in the jump-diffusion
case in \cite{kunitalevy} in the case of first order differentiability,
however, the extension to higher orders of differentiability is not trivial
and requires us to use the framework of \cite{bichteler}. 
\\ 

\noindent
Our proofs rely on the Burkholder-Davis-Gundy (BDG) 
inequality, which states that for any martingale $(M_t)_{t\in \real_+}$
and for all $p \geq 1$ we have 
\begin{equation}
\label{eq:BDG}
\ee \big[ \big|M_t^*\big|^p \big]
\leq C_p \ee\big[[M,M]^{p/2}_t \big], \qquad t\in \real_+, 
\end{equation} 
where $M_t^*=\sup_{s\in [0,t]}|M_s|$, with 
\begin{equation}
  \label{cp}
  C_p = (10 p)^{p/2} \mbox{ ~for } 1\leq p < 2,
  \ C_2 = 2^p, \mbox{ and } C_p = p^p \left( \frac{e}{2} \right)^{p/2}
  \mbox{ ~for } p > 2,
  \end{equation} 
cf. e.g. Theorem 4.2.12 of \cite{bichtelerbk}
or Theorem~48 in Chapter~IV of \cite{protterb2005}. 
 When $p=2$ we have 
\begin{equation} 
  \label{dklsdf}
  \ee\bigg[ \sup_{s\in [0,t]}|M_s|^p \bigg]
\leq C_p \ee\big[\langle M,M\rangle^{p/2}_t \big], \qquad t\in \real_+, 
\end{equation} 
 which implies the bound 
$$ 
\ee\bigg[ \sup_{s\in [0,t]}|M_s|^p \bigg]
\leq 
C_p \big(
\ee [ 
 \langle M, M \rangle_t 
] 
\big)^{p/2}
,
\quad t\in \real_+, 
$$ 
for $p\in [1,2]$. However, \eqref{dklsdf} does not extend
to any $p>2$, see e.g. Remark~357 page 384 of \cite{situ}. 
For this reason we use Kunita's BDG inequality for jump processes, see  
Theorem~2.11 in \cite{kunitalevy}, which is
recovered under a form similar to Corollary~2.14 in \cite{hausenblas}, 
see Lemma~\ref{prop:Bass_Cranston} and Corollary~\ref{prop:Bass_Cranston2} below.
This also extends related results obtained in the
case of a standard Poisson process, in \cite{higham}
see Corollary~1 and Lemma~1 therein, 
with application to the numerical solution of 
SDEs.
\\
   
\noindent 
We proceed by deriving 
moment bounds for the solutions of SDEs 
with jumps of the form 
\eqref{eq:edsX} in Theorem~\ref{theo:moment_affine}.
Similar bounds have been obtained in Theorem~3.2 of \cite{kunitalevy},
however, here we work with random $({\cal F}_t)_{t\geq 0}$-adapted
coefficients and
under weaker integrability conditions. 
Other moment bounds for SDEs with jumps have been derived using \eqref{eq:BDG}
in various works, see for example Lemma~1 in
\cite{xinhongzhang} or
Lemma~2.2 in \cite{zhou-zhang}. 
However, those approaches rely on the
incorrect assumption that \eqref{dklsdf}
holds for any $p\geq 1$. 
Nevertheless, \eqref{dklsdf} is valid for $p=2$, and in this case 
it has been used in \cite{zhou-zhang} to derive
bounds on $\ee\big[ \sup_{s\in [0,t]}|M_s|^p \big]$
for $p > 1$, see the proof of Theorem~2.1 therein 
and also \cite{privault-wang} for an application of
Kunita's BDG inequality to SIR population models for any $p>2$. 
\\

\noindent
The proofs of Proposition~\ref{prop:faa} on the existence of
the flow derivatives and of Theorem~\ref{theo:moment_derivees} on
their $L^p$ integrability 
rely on Theorems~6-29 and 6-44 of \cite{bichteler}.
For this reason, in Sections~\ref{sec:flow_derivative} and
\ref{sec:regularity} we will assume that
the compensators $\nu_t(dy)$, $t\in\real_+ $, in \eqref{eq:edsX} are
dominated by a (deterministic) measure $\eta$ on $\real$, i.e. 
\begin{equation}
\label{eq:compensator} 
\nu_t(A)\leq \eta(A),\quad A\in{\cal B}(\rit), 
\end{equation}
\textcolor{black}{where ${\cal B}(\rit)$ denotes the Borel $\sigma$-algebra of $\rit$},
in addition to the following Assumption (\hyperlink{BGJhyp}{$A_n$}), 
see $(A'\mbox{-}r)$ page~60 \textcolor{black}{in \cite{bichteler}}. 
  
\medskip\noindent
{\bf Assumption (\hypertarget{BGJhyp}{$A_n$}):} 
{\em
 For every $t\in \real_+ $, the functions $r(t,\cdot ): \real\to\real$, $\sigma (t,\cdot ): \real\to\real$
and $g(t,\cdot ) : \real\times\real\to\real$ are ${\cal C}^n$-differentiable and there
is a constant $C > 0$ such that 
$$
\left|\frac{\partial^k r}{\partial x^k}(t,x) \right|
\leq C,
\quad 
\left|\frac{\partial^k \sigma}{\partial x^k}(t,x) \right|
\leq C,
\quad 
 \left|\frac{\partial^{k+l} g}{\partial x^k\partial y^l}(t,x,y) \right| 
 \leq C, 
$$ 
 for all $k,l = 1, \ldots , n$ with $1\leq k+l \leq n$, 
 and a function $\theta\in \bigcap_{q\geq 2} L^q(\real, \eta)$ such that
\begin{equation}
  \label{5.1} 
 \left| \frac{\partial^k g}{\partial x^k}(t,x,y) \right|
 \leq C \theta(y), \qquad x, y\in \real,
 \quad
 k = 1, \ldots , n.
\end{equation} 
}
\vskip-0.2cm 
\noindent
Although the results of \cite{bichteler} are stated for
time-homogeneous SDE coefficients in \eqref{eq:edsX},
they remain valid under our time-inhomogeneous Assumption
(\hyperlink{BGJhyp}{$A_n$}).
For this, we note that Theorems~6-20, 6-24, 6-29 and
6-44 in \cite{bichteler} all rely on Lemma~5.1 page~44 therein, which extends to
the time-inhomogeneous case 
thanks to the domination condition \eqref{eq:compensator},
see Lemma~A.14 of \cite{Bichteler-Jacod83} and
Theorem~2.1 in \cite{Bichteler81}.
\\
  
\noindent
Under \eqref{eq:compensator} and Assumption (\hyperlink{BGJhyp}{$A_n$}), 
in Theorem~\ref{theo:moment_derivees}
we provide sufficient conditions for the flow derivative
\begin{equation}
  \label{fd} 
X^{(n)}_t(x):=\frac{\partial^n X_t}{\partial x^n}(x)
\end{equation} 
to exist and belong to $L^p(\Omega )$
uniformly in $(x,t)\in \real\times [0,T]$,
i.e., 
\begin{equation}
\nonumber 
\sup_{x \in \real } \ee\bigg[
  \sup_{ t \in [0,T]}
  \big|X_t^{(n)}(x)\big|^p\bigg] < \infty, 
\end{equation}
for given orders of derivation $n\geq 1$ and of integrability $p\geq 1$.
\\

\noindent
Flow regularity results up to order four
of differentiability have also been obtained in \cite{prottertalay}
based on a different version of the BDG inequality 
for L\'evy processes with stochastic integrands depending
only on time (see Lemma~4.1 \textcolor{black}{page~409 of \cite{prottertalay}}), 
with application to the convergence of the Euler method. 
\\

\noindent 
\noindent  
As a consequence of \eqref{fd},
when $f(x) : = (x-K)^+$ is the (Lipschitz) payoff function
of a European call option, we can also express the Delta, 
or first derivative of the option price with respect to the
underlying price $x$ as 
$$
\frac{\partial P_t f}{\partial x} (x) = 
\ee\left[
 {\bf 1}_{\{ X_t(x) \geq K \} } X^{(1)}_t(x)
\right]
,
\quad t\in \real_+, \quad x\in \real. 
$$
More generally, 
given the transition semigroup $(P_t)_{t\geq 0}$ of $(X_t(x))_{t\geq 0}$, defined as 
$$
P_t f (x) : = \ee\big[f(X_t(x))\big],
\qquad
t\in \real_+ ,
\quad
x\in \real, 
$$
we deduce 
that for any $f\in {\cal C}_b^\infty (\real)$ the function
$x\mapsto P_tf(x)$ is ${\cal C}^\infty$, with
$$
\frac{\partial^n P_t f}{\partial x^n} (x) = 
\sum_{\pi\in\Pi [n]}
\ee\left[
  \frac{\partial^{|\pi|} f}{\partial x^{|\pi|}}\big(x, X_t(x)\big) \prod_{B\in\pi}X^{(|B|)}_t(x)
\right]
,
\quad t\in \real_+, \quad x\in \real, 
$$
by the Fa\`a di Bruno formula,
where 
the sum over
$\pi$ runs in the set $\Pi [n]$ of all partitions of $\{1, 2, \dots,n\}$,
the product over $B\in\pi$ runs in
all blocks $B$ in the partition $\pi$, 
and $|A|$ stands for the \textcolor{black}{cardinality} of the set $A$. 
 The moment bounds obtained in this paper are also applied to the derivation of
distance estimates between jump-diffusion processes
in \cite{bretonprivault_wasserstein}. 
\\ 

\noindent 
We proceed as follows.
In Section~\ref{sec:BDGjump} we derive
two versions of the BDG inequality with jumps,
similarly to Theorem~2.11 \cite{kunitalevy}
and to Corollary~2.14 of \cite{hausenblas}, and we show that they can be
unified in Corollary~\ref{prop:Bass_Cranston2}. 
This is followed by moment bounds for SDEs in Section~\ref{section:moments}.
In Section~\ref{sec:flow_derivative}
we deal with the flow derivatives $X^{(n)}_t(x)$ by noting that they
satisfy an affine SDE, for which moment bounds can be obtained
from Theorem~\ref{theo:moment_affine},
see Proposition~\ref{prop:faa}. 
Next, in Section~\ref{sec:regularity}
we present our result on moment bounds for flow derivatives,
see Theorem~\ref{theo:moment_derivees}.

\section{Burkholder-Davis-Gundy inequality with jumps}
\label{sec:BDGjump}

\noindent
 Our moment bounds rely on a version of the BDG inequality \eqref{eq:BDG}
which uses the compensator $\langle M, M\rangle $
of $M$ instead of its bracket $[M,M]$.
Consider the compensated Poisson stochastic integral process
\begin{equation}
  \label{Kt} 
K_t:=\int_0^t \int_{-\infty}^\infty g_s (y)\ \big(N(ds,dy)-\nu_s(dy)ds\big),
\qquad t\in \real_+, 
\end{equation} 
of the predictable integrand 
$(g_s(y))_{(s,y)\in \real_+ \times \real}$, 
where $N$ is a Poisson random
measure on $\real_+\times\real$ with compensator $\nu_t(dy)dt$.
When $p=2$, the BDG inequality \eqref{eq:BDG} shows that 
$$ 
\ee\big[\big(K_t^*\big)^2\big] \leq 
2 
\ee\left[\int_0^t\int_{-\infty}^\infty \big( g_s(y) \big)^2 \ N(ds,dy)\right]
=
2 \ee \left[
  \int_0^t\int_{-\infty}^\infty \big(g_s(y)\big)^2 \nu_s(dy)ds
  \right],
$$ 
where $K_t^*=\sup_{s\in [0,t]}|K_s|$, 
since $t\mapsto \int_0^t\int_{-\infty}^\infty \big( g_s(y) \big)^2 \ \big(N(ds,dy)-\nu_s(dy)ds\big)$ is a martingale.  
 In particular, for any $p\in [1,2]$ we have 
 \begin{eqnarray} 
\nonumber 
\ee\big[(K_t^*)^p \big]
& \leq &  
\big(
\ee\big[(K_t^*)^2\big]
\big)^{p/2}
\\
\label{fdfs}
& \leq &  
2^{p/2} \left(
\ee\left[ 
  \int_0^t\int_{-\infty}^\infty (g_s(y))^2 \ \nu_s(dy)ds
  \right]
\right)^{p/2}
,
\quad t\in \real_+ . 
\end{eqnarray} 
Lemma~\ref{prop:Bass_Cranston} below extends the BDG inequality to
$p>2$ with explicit bounding constants,
in relation to the BDG inequality stated for $1<p \leq 2$ in
Corollary~2.14 of \cite{hausenblas}. 
\begin{lemma} 
\label{prop:Bass_Cranston}
Consider the compensated Poisson stochastic integral process 
$(K_t)_{t\in \real_+}$ in \eqref{Kt} 
of a predictable integrand 
$(g_s(y))_{(s,y)\in \real_+ \times \real}$. 
Then, for all $p \geq 2$ we have 
\begin{eqnarray} 
\label{eq:BDGPoisson}
\ee\big[(K_t^*)^p \big]
& \leq &
\frac{2}{p}
(40p)^{p/2}
\left( \frac{p^2 e}{2}\right)^{p ( \log_2p ) /2}
 \ee\left[\int_0^t\int_{-\infty}^\infty | g_s(y) |^p \ \nu_s (dy) ds \right]
 \\
 \nonumber 
 & & + 
 2^p 
 \sum_{k=1}^{\lceil \log_2p\rceil -1}
 \!
\frac{p^{pk}}{2^k} 
 \left(\frac{e}{2}\right)^{kp/2}
 \ 
\ee\left[ \Big(
  \int_0^t\int_{-\infty}^\infty (g_s(y))^{2^k}\ \nu_s(dy)ds\Big)^{p/2^k} \right],
\quad t\in \real_+. 
\end{eqnarray} 
\end{lemma} 
\begin{Proof}
For $r \geq 1$ let 
$$
K^{(r)}_t:=\int_0^t\int_{-\infty}^\infty (g_s(y))^r \ \big(N(ds,dy)-\nu_s(dy)ds\big),
$$
with $K^{(1)}_t=K_t$, $t\in \real_+$. 
 When $p>2$, since $x\mapsto |x|^{p/2}$ is convex, \eqref{eq:BDG} entails
\begin{align} 
\nonumber
 & 
  \ee\big[\big(\big(K_t^{(r)}\big)^*\big)^p\big]\leq
p^p \left(\frac{e}{2}\right)^{p/2}  \ee\left[\Big(\int_0^t\int_{-\infty}^\infty \big( g_s(y) \big)^{2r}\ N(ds,dy)\Big)^{p/2}\right]  
\\
\nonumber
& = p^p \left( \frac{e}{2} \right)^{p/2}
\ee\left[\left|\int_0^t\int_{-\infty}^\infty \big( g_s(y) \big)^{2r}\ \big(N(ds,dy)-\nu_s(dy)\big)ds+ \int_0^t\int_{-\infty}^\infty \big( g_s(y) \big)^{2r}\ \nu_s(dy)ds\right|^{p/2}\right] \\
\nonumber
&\leq  
2^{p/2-1} p^p \left( \frac{e}{2}\right)^{p/2}    
\ee\left[\left|\int_0^t\int_{-\infty}^\infty \big( g_s(y) \big)^{2r}\ \big(N(ds,dy)-\nu_s(dy) ds \big)\right|^{p/2}\right]
\\
\nonumber
& + 2^{p/2-1} \left( \frac{e}{2}\right)^{p/2} 
\ee\left[ \left|\int_0^t\int_{-\infty}^\infty \big( g_s(y) \big)^{2r}\ \nu_s(dy)ds\right|^{p/2}\right] 
\\
\nonumber
&\leq 2^{p/2-1} p^p
\left( \frac{e}{2}\right)^{p/2}
 \ee\big[\big(\big(K_t^{(2r)}\big)^*\big)^{p/2} \big]
 +2^{p/2 -1} p^p
 \left( \frac{e}{2}\right)^{p/2}
 \ee\left[ \Big(\int_0^t\int_{-\infty}^\infty \big( g_s(y) \big)^{2r}\ \nu_s(dy)ds\Big)^{p/2}\right].
 \\
 \label{eq:BDGrec}
\end{align} 
The recursive application of the bound \eqref{eq:BDGrec} starting from
$K^{(1)}_t=K_t$ yields 
\begin{eqnarray*} 
\ee\big[(K_t^*)^{p}\big]
& \leq &  
 p^{np} \left( \frac{e}{2}\right)^{pn/2} \bigg(\prod_{j=1}^n 2^{p/2^j-1}\bigg)\
 \ee\big[\big(\big(K_t^{(2^n)}\big)^*\big)^{p/2^n}\big]
 \\
  & & + 
 \sum_{k=1}^n p^{kp}
 \left( \frac{e}{2}\right)^{pn/2}
 \bigg(\prod_{j=1}^k 2^{p/2^j-1}\bigg)\
\ee\left[
  \Big(
  \int_0^t\int_{-\infty}^\infty (g_s(y))^{2^k}\ \nu_s(dy)ds\Big)^{p/2^k}
  \right].
\end{eqnarray*}
Taking $n = \lceil \log_2p\rceil -1$,
i.e. $p \in (2^n,2^{n+1}]$, by \eqref{eq:BDG} we have 
 \begin{eqnarray*} 
 \ee\big[\big(\big(K_t^{(2^n)}\big)^*\big)^{p/2^n}\big]
& \leq & 
 (10p)^{p/2} 
 \ee\left[\Big(\int_0^t\int_{-\infty}^\infty \big( g_s(y) \big)^{2^{n+1}}\ N(ds,dy)\Big)^{p/2^{n+1}}\right]  
 \\
 & \leq & 
 (10p)^{p/2}
 \ee\left[\int_0^t\int_{-\infty}^\infty | g_s(y) |^p \ N(ds,dy)\right]
  \\
 & = & 
  (10p)^{p/2}
  \ee\left[\int_0^t\int_{-\infty}^\infty | g_s(y) |^p \ \nu_s(dy) ds \right],  
\end{eqnarray*} 
 since $t\mapsto \int_0^t\int_{-\infty}^\infty \big( g_s(y) \big)^2 \ \big(N(ds,dy)-\nu_s(dy)ds\big)$ is a martingale, 
where we used the fact that
$$
\textcolor{black}{
  \Vert a \Vert_{\ell^2}:
= \left(
\sum_{k=0}^\infty
(a_k)^2 \right)^{1/2}
\leq
 \left(
\sum_{k=0}^\infty
(a_k)^q \right)^{1/q}
= : \Vert a \Vert_{\ell^q}
}
$$
 for any real sequence $(a_k)_{k\in \inte}$ and $1 \leq q = p/2^n \leq 2$
 as on page~410 after Equation~(22) in \cite{prottertalay} ,
 which allows us to conclude to \eqref{eq:BDGPoisson}.
\end{Proof} 
 From Lemma~\ref{prop:Bass_Cranston} we recover the
 following version of the Kunita's BDG inequality with jumps,
 cf. Theorem~2.11 of \cite{kunitalevy} and
 Theorem~4.4.23 of \cite{applebk2}. 
\begin{corollary} 
\label{prop:Bass_Cranston2}
Consider the compensated Poisson stochastic integral process 
\begin{equation} 
K_t:=
x + \int_0^tu_s\ ds + 
\int_0^tv_s\ dW_s + 
\int_0^t \int_{-\infty}^\infty g_s (y)\ \big(N(ds,dy)-\nu_s(dy)ds\big),
\quad t\in \real_+, 
\\
\label{djkldsf}
\end{equation} 
of the predictable integrands 
$(u_s)_{s\in \real_+}$,
$(v_s)_{s\in \real_+}$,
$(g_s(y))_{(s,y)\in \real_+ \times \real}$. 
Then, for all $p \geq 2$ and
$T\in \real_+$ we have 
\begin{align*} 
\nonumber 
 & \ee\big[(K_T^*)^p \big]
   \\
   & \leq 
   2^{2p-2}
   \left(
   |x|^p + 
   \ee\left[ \left( \int_0^T |u_t| \ dt \right)^p \right]
  +
   C_p
  \ee\left[\Big(\int_0^T |v_t|^2\ dt\Big)^{p/2}\right] 
  \right.
  \\
  & \quad 
  \left.
  +
   \widetilde{C}_p
\ee\left[ 
  \int_0^T\int_{-\infty}^\infty |g_t(y)|^p \ \nu_t(dy)dt
  \right]
+
\widetilde{C}_p
\ee\left[ \left(
  \int_0^T\int_{-\infty}^\infty |g_t(y)|^2 \ \nu_t(dy)dt\right)^{ p/2 } 
  \right]
\right) 
,
\end{align*} 
 where $C_p$ is defined in \eqref{cp}, and 
$$\widetilde{C}_p = 
 \frac{2}{p}
(40p)^{p/2}
\left( \frac{p^2 e}{2}\right)^{p ( \log_2p ) /2}
  + 
 2^p 
 \sum_{k=1}^{\lceil \log_2p\rceil -1}
\frac{p^{pk}}{2^k} 
 \left(\frac{e}{2}\right)^{kp/2}
\leq 
2^p
p^{ p \log_2p}
\left(
2 + \left( 10 e^{\lceil \log_2p \rceil} \right)^{p/2}  
\right). 
$$
\end{corollary} 
\begin{Proof}
 By the convexity of $\real \ni x \mapsto |x|^p$,  we have 
\begin{eqnarray}
\nonumber
  \ee\left[\sup_{s\in [0,T]} |K_{t^-}|^p\right] 
 & \leq & 4^{p-1}\bigg(x^p
+\ee\left[\sup_{t\in [0,T]} \Big|\int_0^tu_s\ ds\Big|^p\right]
+\ee\left[\sup_{t\in [0,T]} \Big|\int_0^tv_s\ dW_s\Big|^p\right]
\\
\nonumber
& & 
\left.
+\ee\bigg[ \sup_{t\in [0,T]} \Big|\int_0^{t^-}\int_{-\infty}^\infty g_{s^-}(z)\ \big(N(ds,dz)-\nu_s(dz)ds\big)\Big|^p\bigg] \right). 
\end{eqnarray} 
Further, by the (standard) BDG inequality
for Brownian stochastic integrals, 
 the Jensen inequality for the uniform measure on $[0,t]$
and the Fubini theorem, we find: 
$$ 
  \ee\left[\sup_{t\in [0,T]} \Big|\int_0^tv_s\ dW_s\Big|^p\right]
  \leq C_p \ee\left[\Big(\int_0^T |v_s|^2\ ds\Big)^{p/2}\right].
$$
  Regarding the jump term, by the log-convexity in $p$ of the $L^p$ norms,
  taking $n = \lceil \log_2p\rceil -1$, 
  i.e. $p \in (2^n,2^{n+1}]$,
  and $\theta \in (0,1)$ such that
$2^{-k}= (1-\theta) / p + \theta / 2$,
    we have 
  \begin{eqnarray*} 
  \lefteqn{
    \ee\left[ \Big(
  \int_0^T\int_{-\infty}^\infty (g_t(y))^{2^k}\ \nu_t(dy)dt\Big)^{p/2^k}\right]
  }
  \\
  & \leq & 
\ee\left[ \left(
  \int_0^T\int_{-\infty}^\infty (g_t(y))^p \ \nu_t(dy)dt\right)^{1-\theta}
  \left(
  \int_0^T\int_{-\infty}^\infty (g_t(y))^2 \ \nu_t(dy)dt\right)^{\theta p/2 } 
  \right]
\\
& \leq & 
\left(
\ee\left[ 
  \int_0^T\int_{-\infty}^\infty (g_t(y))^p \ \nu_t(dy)dt
  \right]
\right)^{1-\theta}
\left(
\ee\left[ \left(
  \int_0^T\int_{-\infty}^\infty (g_t(y))^2 \ \nu_t(dy)dt\right)^{ p/2 } 
  \right]
\right)^{\theta}
\\
& \leq & 
  (1-\theta) 
\ee\left[ 
  \int_0^T\int_{-\infty}^\infty (g_t(y))^p \ \nu_t(dy)dt\
  \right]
+
\theta 
\ee\left[ \left(
  \int_0^T\int_{-\infty}^\infty (g_t(y))^2 \ \nu_t(dy)dt\right)^{ p/2 } 
  \right]
\\
& \leq & 
\ee\left[ 
  \int_0^T\int_{-\infty}^\infty (g_t(y))^p \ \nu_t(dy)dt
  \right]
+
\ee\left[ \left(
  \int_0^T\int_{-\infty}^\infty (g_t(y))^2 \ \nu_t(dy)dt\right)^{ p/2 } 
  \right]
,
\quad T\in \real_+, 
\end{eqnarray*} 
$k=1,\ldots , n$, after using the H\"older inequality and the
bound 
$x^{1-\theta} y^\theta \leq (1-\theta )x + \theta y$,
$x,y \geq 0$.
Hence, substituting this bound in \eqref{eq:BDGPoisson}, we obtain 
 \begin{eqnarray*} 
  \lefteqn{
\ee\bigg[ \sup_{t\in [0,T]} \Big|\int_0^{t^-}\int_{-\infty}^\infty g_{s^-}(z)\ \big(N(ds,dz)-\nu_s(dz)ds\big)\Big|^p\bigg] 
}
 \\
 \nonumber & \leq &  
 \frac{2^{p+1}}{p} p^{p\log_2 p} (10p)^{p/2} \left( \frac{e}{2}\right)^{p ( \log_2p ) /2}
 \ee\left[\int_0^T \int_{-\infty}^\infty \big( g_t(y) \big)^p \ \nu_t (dy) dt \right]
\\
\nonumber 
 & & + 
 2^p 
 \sum_{k=1}^{\lceil \log_2p\rceil -1}
 \frac{p^{pk}}{2^k} 
 \left(\frac{e}{2}\right)^{kp/2}
\ee\left[ \Big(
  \int_0^T \int_{-\infty}^\infty (g_t(y))^{2^k}\ \nu_t(dy)dt\Big)^{p/2^k} \right]
  \\
  & \leq &
  \widetilde{C}_p
  \left(
\ee\left[ 
  \int_0^T\int_{-\infty}^\infty (g_t(y))^p \ \nu_t(dy)dt
  \right]
+
\ee\left[ \left(
  \int_0^T\int_{-\infty}^\infty (g_t(y))^2 \ \nu_t(dy)dt\right)^{ p/2 } 
  \right]
\right). 
\end{eqnarray*} 
\end{Proof} 
The following consequence of Corollary~\ref{prop:Bass_Cranston2}
recovers Corollary~2.12 in \cite{kunitalevy} using the H\"older inequality. 
\begin{corollary} 
\label{prop:Bass_Cranston3}
Consider the compensated Poisson stochastic integral process 
$(K_t)_{t\in [0,T]}$ in \eqref{djkldsf}. For all $p \geq 2$ and
$T\in \real_+$ we have 
\begin{align*} 
\nonumber 
 & \ee\big[(K_T^*)^p \big]
 \leq 
   2^{2p-2}
\left(   |x|^p + 
  T^{p-1} \ee\left[\int_0^T |u_t|^p \ dt \right]
  +
   C_p T^{p/2-1} 
  \ee\left[ \int_0^T |v_t|^p\ dt \right] 
  \right.
  \\
  & 
  \left. +
  \widetilde{C}_p
\ee\left[ 
  \int_0^T\int_{-\infty}^\infty |g_t(y)|^p \ \nu_t(dy)dt
  \right]
+
 T^{p/2-1} \widetilde{C}_p
\ee\left[ 
  \int_0^T \left( \int_{-\infty}^\infty |g_t(y)|^2 \ \nu_t(dy)\right)^{ p/2 }
  dt
  \right]
\right).
\end{align*} 
\end{corollary} 
When the integrand $g_t(y)$ satisfies
$|g_t(y)| \leq |f(y)| |g_t|$ where $f(y)$
is a deterministic function of $y\in \real$, 
$(g_t)_{t\in \real_+}$ is an $({\cal F}_t)_{t\geq 0}$-adapted process, 
and $(\nu_t (dy))_{t \in \real_+ }=\nu (dy)$, $t\in \real_+$, 
is the intensity measure of a time-homogeneous Poisson point process, 
Corollary~\ref{prop:Bass_Cranston3} yields
\begin{align*} 
\nonumber 
 & \ee\big[(K_T^*)^p \big]
 \leq 
   2^{2p-2}
\left(   |x|^p + 
  T^{p-1} \ee\left[\int_0^T |u_t|^p \ dt \right]
  +
   C_p T^{p/2-1} 
  \ee\left[ \int_0^T |v_t|^p\ dt \right] 
  \right.
  \\
  & 
  \left. +
  \widetilde{C}_p
  \left( 
  \int_{-\infty}^\infty |f(y)|^p \nu (dy)
+
 T^{p/2-1} 
 \left( \int_{-\infty}^\infty |f(y)|^2 \ \nu (dy)\right)^{ p/2 }
 \right)
 \ee\left[ 
   \int_0^T
   |g_t|^p
 dt
  \right]
\right), 
\end{align*} 
which recovers related versions of the BDG inequality
with jumps 
such as Lemma~5.2 of \cite{bass} 
which is stated for $p=2^n$, $n\geq 1$, 
or Lemma~4.1 of \cite{prottertalay} which is stated for $p \geq 2$
using a related recursion. 
We also refer the reader to 
Lemma~A.14 of \cite{Bichteler-Jacod83} and to 
the proof of Theorem~2.1 in \cite{Bichteler81}, 
or to \cite{lepingle2} and \cite{pratelli2}, 
for other versions of the BDG inequality with jumps.
\section{Moment bounds for SDE solutions} 
\label{section:moments}

In this section, we derive moment bounds for jump-diffusion
SDEs, based on the BDG inequality with jumps given in
Corollary~\ref{prop:Bass_Cranston2}. 
\\
 
\noindent
The following result provides moment bounds in $L^p(\Omega )$,
$p \geq 2$, on the solution of SDEs of the form
\begin{equation} 
  dX_t = a_t ( X_t) dt +b_t (X_t ) dW_t
+\int_{-\infty}^\infty c_{t^-} (z , X_{t^-} ) \big(N(dt,dz)-\nu_t(dz)dt\big),
\\
\label{eq:nonaffine}
\end{equation} 
whose existence and uniqueness of solutions
follows by standard arguments, see e.g. Theorem~3.1 in \cite{kunitalevy}. 
In contrast with Theorem~3.2 of \cite{kunitalevy},
we work with random $({\cal F}_t)_{t\geq 0}$-adapted coefficients and
under weaker integrability conditions as 
Condition~(3.2) in \cite{kunitalevy} requires integrability
of all orders. 
We let $\| X \|_\infty$ stand for the
$L^\infty(\Omega)$ norm of a random variable $X$.
\begin{theorem}
\label{theo:moment_affine}
Let $p\geq 2$ and
consider the solution $(X_t)_{t\geq 0}$ of
the one-dimensional solution of the jump-diffusion SDE 
\eqref{eq:nonaffine}, where 
 the coefficients
 $(a_t(x))_{t\in [0,T]}$, $(b_t(x))_{t\in [0,T]}$, $(c_t(z,x))_{t\in [0,T]}$
 are $({\cal F}_t)_{t\geq 0}$-adapted processes such that 
$$
| a_t ( x) - a_t(y) | \leq a_t | x-y|,
\quad 
| b_t( x) - b_t(y) | \leq b_t | x-y|,
 \qquad x,y \in \real, 
$$
 where
$(a_t(0))_{t\in [0,T]} ,
(b_t(0))_{t\in [0,T]}
 \in L^p ( \Omega \times [0,T])$,
 $(a_t)_{t\in [0,T]} ,
(b_t)_{t\in [0,T]}
\in L^p ( [0,T], L^\infty (\Omega ))$,
 and
 $$
 | c_t (z , x ) - c_t ( z , y ) | \leq c_t(z) |x-y|,
 \quad x, y \in \real,
 $$ 
where
\begin{equation} 
\label{eq:cc_hyp}
\int_0^T
\Big\|\int_{-\infty}^\infty
 | c_t (z) |^2 
 \ \nu_t (dz)\Big\|_\infty^{p/2} \ dt < \infty, 
 \ \ 
 \int_0^T
\Big\|\int_{-\infty}^\infty
| c_t (z) |^p 
\ \nu_t (dz)\Big\|_\infty \ dt < \infty, 
\end{equation} 
and
\begin{equation} 
\label{eq:cc_hyp2}
\ee\left[ \int_0^T
  \Big( \int_{-\infty}^\infty
  | c_t(z,0) |^2 \ \nu_t(dz)
  \Big)^{p/2}
 \ dt\right]
< \infty,
\ 
\ee\left[ \int_0^T
 \int_{-\infty}^\infty
 | c_t(z,0) |^p 
 \ \nu_t(dz)  \ dt\right] < \infty. 
\end{equation} 
 Then we have 
$$
 \ee\bigg[ \sup_{t\in[0,T]} |X_t|^p\bigg]<C(p,T) <\infty, 
$$
where $C(p,T)$ depends on the above norms of $a, b, c, u, v, w$. 
\end{theorem}
\begin{Proof} 
 We have  
\begin{equation} 
\label{eq:dd2}
\E \left[
  \int_0^T | a_t(X_t) |^p dt \right]
\leq 2^{p-1} \left( \int_0^T \|a_t\|_\infty^p \ee\big[|X_t|^p\big] dt
+ 
\E \left[ \int_0^T |a_t(0)|^p dt \right] 
\right) ,
\end{equation} 
and
\begin{equation} 
\label{eq:mb2}
\E \left[
  \int_0^T | b_t(X_t) |^p dt \right]
\leq 2^{p-1} \left( \int_0^T \|b_t\|_\infty^p \ee\big[|X_t|^p\big] dt
+ 
\E \left[ \int_0^T |b_t(0)|^p dt \right] 
\right). 
\end{equation} 
Regarding the jump term, we note that 
\begin{align} 
\nonumber
& 
 \widetilde{C}_p 
 T^{p/2-1}
 \ee\left[ \int_0^T
   \Big(
   \int_{-\infty}^\infty
    | c_t(z,X_t) |^2
  \ \nu_t(dz) \Big)^{p/2} \ dt\right]
+ 
 \widetilde{C}_p \ee\left[ \int_0^T
  \int_{-\infty}^\infty
    | c_t(z, X_t) |^p
  \ \nu_t(dz) \ dt\right]
\\
\nonumber
&\leq  
2^p 
 \widetilde{C}_p
T^{p/2-1}
\ee\left[ \int_0^T\Big(
  \int_{-\infty}^\infty
  | c_t(z) |^2 |X_t|^2 \ \nu_t(dz)
  +
  \int_{-\infty}^\infty
  | c_t(z,0) |^2 \ \nu_t(dz)\Big)^{p/2} \ dt\right]
\\
\nonumber
& \quad 
+ 2^{p-1} 
 \widetilde{C}_p
\ee\left[ \int_0^T
 \int_{-\infty}^\infty
 \big( | c_t(z) |^p |X_t|^p + | c_t(z,0) |^p \big)
 \ \nu_t(dz)  \ dt\right]
\\
\nonumber
&\leq  
2^{3p/2-1} 
 \widetilde{C}_p
T^{p/2-1}
\ee\left[ \int_0^T
  \Big(
  \Big(  \int_{-\infty}^\infty
  | c_t(z) |^2 \ \nu_t(dz)
  \Big)^{p/2} |X_t|^p 
  +
  \left( \int_{-\infty}^\infty
  | c_t(z,0) |^2 \ \nu_t(dz)
  \Big)^{p/2}
  \right) \ dt\right]
\\
\nonumber
& \quad 
+
2^{p-1}
\widetilde{C}_p
\ee\left[ \int_0^T
 \int_{-\infty}^\infty
 | c_t(z) |^p \ \nu_t(dz) |X_t |^p \ dt\right]
+
2^{p-1} \widetilde{C}_p
\ee\left[ \int_0^T
 \int_{-\infty}^\infty
 | c_t(z,0) |^p 
 \ \nu_t(dz)  \ dt\right]
\\
\label{eq:NN2}
&\leq  
2^{3p/2-1} 
 \widetilde{C}_p
T^{p/2-1}\int_0^T\left\|\int_{-\infty}^\infty
 | c_t(z) |^2 
\ \nu_t(dz)\right\|_\infty^{p/2}\ee\big[ X_t^p \big]\ dt
\\
\nonumber
& \quad  +
2^{3p/2-1} 
 \widetilde{C}_p
T^{p/2-1}
\ee\left[ \int_0^T
  \Big( \int_{-\infty}^\infty
  | c_t(z,0) |^2 \ \nu_t(dz)
  \Big)^{p/2}
 \ dt\right]
\\
\nonumber
 & \quad  
+ 2^{p-1} 
 \widetilde{C}_p
\int_0^T\left\|\int_{-\infty}^\infty
 | c_t(z) |^p 
\ \nu_t(dz)\right\|_\infty \ee\big[ X_t^p\big]\ dt
+
2^{p-1} \widetilde{C}_p
\ee\left[ \int_0^T
 \int_{-\infty}^\infty
 | c_t(z,0) |^p 
 \ \nu_t(dz)  \ dt\right].
\end{align} 
\noindent
Hence by the BDG inequality of Corollary~\ref{prop:Bass_Cranston3}
and the bounds 
\eqref{eq:dd2}, \eqref{eq:mb2}, 
 \eqref{eq:NN2}, setting
\begin{align*}
 & 
 F(T):= 4^{p-1}\bigg(
x^p 
  + (2T)^{p-1} \E \left[ \int_0^T |a_t(0)|^p dt \right]  
+ 2^{p-1} C_p T^{p/2-1} \E \left[ \int_0^T |b_t(0)|^p dt \right]  
 \\
& 
\left.
+
2^{p-1} 
 \widetilde{C}_p
 \left(
 2^{p/2}
 T^{p/2-1}
\ee\left[ \int_0^T
  \Big( \int_{-\infty}^\infty
  | c_t(z,0) |^2 \ \nu_t(dz)
  \Big)^{p/2}
 \ dt\right]
+ 
\ee\left[ \int_0^T
 \int_{-\infty}^\infty
 | c_t(z,0) |^p 
 \ \nu_t(dz)  \ dt\right]
\right)
\right),
\end{align*} 
and 
\begin{align*}
 & 
    G(t) : = 4^{p-1}\Big(
(2T)^{p-1} \|a_t\|_\infty^p 
+ 2^{p-1} C_p T^{p/2-1} \|b_t\|_\infty^p 
  \\
& 
\quad +
2^{3p/2-1}
\widetilde{C}_p 
T^{p/2-1}\left\|\int_{-\infty}^\infty
 | c_t(z) |^2 
 \ \nu_t(dz)\right\|_\infty^{p/2} 
\left.
 + 2^{p-1} 
 \widetilde{C}_p 
\left\|\int_{-\infty}^\infty
 | c_t(z) |^p 
 \ \nu_t(dz)\right\|_\infty
 \right),
\end{align*} 
$t\in [0,T]$, we have 
\begin{eqnarray}
\ee\left[\sup_{t\in [0,T]} |X_{t^-}|^p\right]
\nonumber 
&\leq& F(T)+\ee \left[ \int_0^TG(t) |X_t|^p \ dt \right] \\
\nonumber 
&\leq& F(T)+\int_0^TG(t)\ee\left[\sup_{s\in [0,t]}|X_s|^p\right]\ dt, 
\end{eqnarray}
hence by 
the Gr\"onwall lemma we find 
\begin{equation}
\label{eq:momentZ-}
\ee\left[\sup_{t\in [0,T]} |X_{t^-}|^p\right]
\leq C(p,T): = F(T)\exp\left(\int_0^T G(t)\ dt\right), 
\end{equation} 
which is finite 
since $(a_t)_{t\in [0,T]} , (b_t)_{t\in [0,T]} \in L^p ( [0,T], L^\infty (\Omega ))$
and
\eqref{eq:cc_hyp}-\eqref{eq:cc_hyp2} 
are in force. 
Since $X_{t^-}=X_t$ almost surely, the same bound follows for the moment of order $p$ of $X_t$. 
\end{Proof} 
\noindent
Theorem~\ref{theo:moment_affine} applies, in particular, to the solution
$(X_t)_{t\geq 0}$ of the one-dimensional jump-diffusion affine SDE 
\begin{eqnarray} 
  \label{eq:affine}
  dX_t & = & u_t dt+a_t X_t dt
+v_t dW_t +b_t X_tdW_t
\\
\nonumber
& & +\int_{-\infty}^\infty w_{t^-} (z) \big(N(dt,dz)-\nu_t(dz)dt\big)
+X_{t^-} \int_{-\infty}^\infty c_{t^-} (z) \big(N(dt,dz)-\nu_t(dz)dt\big),
\end{eqnarray} 
 by taking
\begin{equation}
  \label{vjdsgf}
  a_t(y)=a_ty+u_t, \quad b_t(y)=b_ty+v_t, \quad c_t(z,y)=c_t(z)y+w_t(z),
\end{equation} 
where $u$, $v$ are in a certain $L^p$ space and $a$, $b$, $c$ are in $L^\infty$.  
\\
 
\noindent
The following uniform version of Theorem~\ref{theo:moment_affine}
will be required in Section~\ref{sec:flow_derivative}.
When the processes 
$a_\alpha$, $b_\alpha$, $c_\alpha$ all depend on a parameter $\alpha\in A$, the solution
$(X_{\alpha ,t})_{t\geq 0}$ of the corresponding SDE \eqref{eq:affine-alpha} below
 enjoys the following uniform bound. 
\begin{corollary}
\label{corol:moment_affineU}
Let $p \geq 2$. 
Assume that 
 the coefficients
 $(a_{\alpha , t} (x))_{t\in [0,T]}$, $(b_{\alpha , t} (x))_{t\in [0,T]}$, $(c_{\alpha , t} (z,x))_{t\in [0,T]}$
 are $({\cal F}_t)_{t\geq 0}$-adapted processes such that 
$$
| a_{\alpha , t} ( x) - a_{\alpha , t} (y) | \leq a_{\alpha , t} | x-y|,
\quad 
| b_{\alpha , t} (x) - b_{\alpha , t} (y) | \leq b_{\alpha , t} | x-y|,
 \qquad
x,y \in \real, 
$$
 where
 $(a_{\alpha , t} (0))_{t\in [0,T]}, (b_{\alpha , t} (0))_{t\in [0,T]}
 \in L^p ( \Omega \times [0,T] )$, 
$(a_{\alpha , t} )_{t\in [0,T]}$, 
$(b_{\alpha , t})_{t\in [0,T]}
\in L^p ( [0,T], L^\infty (\Omega) )$,
uniformly in $\alpha \in A$, with 
$$
 | c_{\alpha , t} (z , x ) - c_{\alpha , t} ( z , y ) | \leq c_{\alpha , t} (z) |x-y|,
 \quad x, y \in \real,
 $$ 
 with 
\begin{equation} 
\label{eq:cc_hypU}
 \sup_{\alpha\in A}\int_0^T\Big\|\int_{-\infty}^\infty
\big|c_{\alpha , t}(z)\big|^2
\ \nu_t(dz)\Big\|_\infty^{p/2}\ dt<\infty,
\ \ 
\sup_{\alpha\in A}\int_0^T\Big\|\int_{-\infty}^\infty
\big|c_{\alpha , t}(z)\big|^p
\ \nu_t(dz)\Big\|_\infty \ dt<\infty, 
\end{equation} 
and
\begin{equation} 
\label{eq:cc_hypU2}
 \sup_{\alpha\in A}
 \E \left[\int_0^T\Big( 
   \int_{-\infty}^\infty
\big|c_{\alpha , t}(z,0)\big|^2
\ \nu_t(dz) \Big)^{p/2}\ dt
\right]
 <\infty,  
\ \ 
\sup_{\alpha\in A}
\E \left[\int_0^T
  \int_{-\infty}^\infty
\big|c_{\alpha , t}(z,0)\big|^p
\ \nu_t(dz) \ dt \right] <\infty. 
\end{equation} 
 Then, for the solutions $X_\alpha$ of the SDE
\begin{equation} 
\label{eq:affine-alpha}
dX_{\alpha , t}  =  
a_{\alpha , t} (X_{\alpha , t}) dt
+b_{\alpha , t} (X_{\alpha , t} ) dW_t
+ \int_{-\infty}^\infty c_{\alpha , t^-} (z, X_{\alpha , t^-}) \big(N(dt,dz)-\nu_t(dz)dt\big),
\end{equation} 
 $\alpha\in A$, we have 
\begin{equation}
\label{eq:bound_affineU}
\sup_{\alpha\in A } \ee\bigg[
  \sup_{t\in[0,T]}
  |X_{\alpha , t}|^p\bigg] < C(p,T) <\infty, 
\end{equation}
where $C(p,T)$ depends on the above norms of $a_\alpha$,
$b_\alpha$, $c_\alpha$,
which are all assumed to be bounded
uniformly in $\alpha\in A$. 
\end{corollary}
\begin{Proof}
  Only the conclusion of the previous proof for Theorem~\ref{theo:moment_affine}
  is required to be changed. 
The bound \eqref{eq:momentZ-} still holds true
for $X_{\alpha , t}$ with the functions
\begin{align*}
& F_\alpha ( T):= 4^{p-1}\bigg(
x^p 
+
(2T)^{p-1} \E \left[ \int_0^T |a_{\alpha , t}(0)|^p dt \right]  
+ 2^{p-1} C_p T^{p/2-1} \E \left[ \int_0^T | b_{\alpha , t}(0) |^p dt \right] 
 \\
\nonumber
&  
\left.
 + 2^{p-1} \widetilde{C}_p\
 \left(
 2^{p/2}
 T^{p/2-1}
 \E \left[ \int_0^T \Big( \int_{-\infty}^\infty
| c_{\alpha , t}(z,0) |^2
\ \nu_t(dz) \Big)^{p/2} dt \right] 
+ 
\E \left[ \int_{-\infty}^\infty
| c_{\alpha , t}(z,0) |^p
\ \nu_t(dz) dt \right] 
\right)
\right)
\end{align*}
and
\begin{align*}
 &      G_\alpha(t)  := 
 4^{p-1}\bigg(
(2T)^{p-1} \|a_{\alpha , t}\|_\infty^p
+2^{p-1} C_p T^{p/2-1} \|b_{\alpha , t}\|_\infty^p
\\
 & 
\qquad + 2^{3p/2-1} \widetilde{C}_p\
T^{p/2-1}\Big\|\int_{-\infty}^\infty
| c_{\alpha , t}(z) |^2
\ \nu_t(dz)\Big\|_\infty^{p/2} 
+ 2^{p-1} \widetilde{C}_p\
\Big\|\int_{-\infty}^\infty
| c_{\alpha , t}(z) |^p
\big)
\ \nu_t(dz)\Big\|_\infty
\bigg), 
\end{align*}
$t\in [0,T]$.
Under the conditions of Corollary~\ref{corol:moment_affineU}, we have 
$$
\sup_{\alpha\in A} F_\alpha (T)<\infty, 
\quad 
\sup_{\alpha\in A} G_\alpha (t)<\infty,\qquad t\in [0,T],
$$
and the conclusion \eqref{eq:bound_affineU} follows likewise.
\end{Proof}

\section{Flow derivatives} 
\label{sec:flow_derivative}

In this section, we show that
the derivatives of the flow $x\mapsto X_t(x)$ of the SDE
\eqref{eq:edsX} are solutions of an affine SDE. 

\vskip0.2cm 

\noindent
{\bf Convention.} 
{\em Given the gradient 
$\nabla_z=\big(\frac{\partial}{\partial z_1}, \dots, \frac{\partial}{\partial z_d}\big)$
and $F:\real^d\to \real^p$,
we denote
$\nabla_zF(z)=\Big(\frac{\partial F_i}{\partial z_j}(z)\Big)_{\begin{subarray}{l}1\leq i\leq p\\ 1\leq j\leq d\end{subarray}}\in M_{p,d}(\real)$,
and under the identification $M_{p,d}(\real)\sim \real^{pd}$, 
we write 
$$
\nabla_z\big(\nabla_z F\big)(x)=\nabla_z^2 F(x)\in M_{pd,d}(\real)\sim \real^{pd^2}.
$$ 
}

\noindent
By successive differentiation of 
$$
r\big(s, X_s(x)\big), \quad
\sigma\big(s, X_s(x)\big),\quad
g\big(s, X_s(x),y\big)
$$ 
with respect to $x$ 
and applying Theorem~6-29 of \cite{bichteler}
recursively under Assumption (\hyperlink{BGJhyp}{$A_n$}),
we obtain the following result. 
\begin{prop}
\label{prop:faa}
Assume that \eqref{eq:compensator} and (\hyperlink{BGJhyp}{$A_n$}) hold
for some $n\geq 1$.
Then the flow $x\mapsto X_t(x)$
of the solution to the real SDE \eqref{eq:edsX} is 
$n$ times differentiable on $\real$ and, for $k = 1, \ldots , n$, 
$X^{(k)}_t(x):= \displaystyle \frac{\partial^kX_t}{\partial x^k}(x)$ is solution of 
\begin{eqnarray}
\nonumber
\lefteqn{
  dX^{(k)}_t(x) =
  \left(\sum_{\pi\in\Pi [k]} \frac{\partial^{|\pi|} r}{\partial x^{|\pi|}}\big(t, X_t(x)\big) \prod_{B\in\pi}X^{(|B|)}_t(x)\right)\ dt
}
\\
\nonumber
&&+\left( \sum_{\pi\in\Pi [k]} \frac{\partial^{|\pi|} \sigma}{\partial x^{|\pi|}}\big(t, X_t(x)\big) \prod_{B\in\pi}X^{(|B|)}_t(x)\right)\ dW_t\\
\label{eq:faa_Xn}
&&+\int_{-\infty}^\infty\left(\sum_{\pi\in\Pi [k]} \frac{\partial^{|\pi|}g}{\partial x^{|\pi|}}\big(t, X_t(x),y\big) \prod_{B\in\pi}X^{(|B|)}_t(x)\right)\ \big(N(dt,dy)-\nu_t(dy)dt\big), \qquad 
\end{eqnarray}
where, again, the sum over
$\pi$ runs in the set $\Pi [k]$ of all partitions of $\{1, 2, \dots,k\}$,
and the product over $B\in\pi$ runs in
all blocks $B$ in the partition $\pi$. 
\end{prop}
\begin{Proof}
For $n=0$, $X^{(0)}_t(x)=X_t(x)$ and \eqref{eq:faa_Xn} reduces to the SDE \eqref{eq:edsX}. 
For $n=1$, \eqref{eq:faa_Xn} is given by Theorem 6-29 in \cite{bichteler}:
\begin{eqnarray}
\nonumber
dX^{(1)}_t(x)&=&
 X^{(1)}_t\textcolor{black}{(x)} \frac{\partial r}{\partial x}(x,X_t(x)) dt
 +
 X^{(1)}_t\textcolor{black}{(x)} \frac{\partial \sigma}{\partial x}(x,X_t(x)) dW_t
\\
 \label{eq:Z1}
 &&
 + 
 X^{(1)}_t(x)
 \int_{-\infty}^\infty \frac{\partial g}{\partial x}(x,X_t(x),y)
  \big(N(dt,dy)-\nu_t(dy)dt\big).
\end{eqnarray}

\noindent
We continue the proof by induction on $n\geq 1$.
We assume that 
$$
Z^{(n-1)}_t (x):=\left(X_t (x), \frac{\partial X_t }{\partial x}(x), \dots, \frac{\partial^{n-1}X_t }{\partial x^{n-1}}(x)\right)^\intercal 
\in \real^n
$$
is solution of the $n$-dimensional SDE 
\begin{eqnarray}
\nonumber
Z^{(n-1)}_t(x)&=&z^{(n-1)}+\int_0^t r^{(n-1)}(s,Z^{(n-1)}_s)\ ds
+\int_0^t \sigma^{(n-1)}(s,Z^{(n-1)}_s)\ dW_s
\\
\label{eq:Zn-1}
&&+\int_0^t\int_{-\infty}^\infty g^{(n-1)}(t,Z^{(n-1)}_{s^-}(x),y)\ \big(N(ds,dy)-\nu_s(dy)ds \big) 
\end{eqnarray}
with
$$z^{(n-1)}_0=x, \quad
z^{(n-1)}_1=1 \mbox{ ~and~ } z^{(n-1)}_k=0, \qquad k\geq 2, 
$$
 and 
 \begin{enumerate}[a)] 
 \item
   $r^{(n-1)} :[0,T]\times\real^n\to\real^n$ is given by 
 \begin{equation}
   \label{eq:rn}
r^{(n-1)}\big(t,z^{(n-1)}\big)
=
\big(r_0^{(n-1)}, \dots, r_{n-1}^{(n-1)}\big)^\intercal
=\big(r^{(n-2)}\big(t,z^{(n-2)}\big),\ r_{n-1}^{(n-1)}(t, z^{(n-1))}\big)^\intercal,
\end{equation} 
 where 
\begin{equation} 
\label{eq:faa_r}
 r^{(n-1)}_{n-1}\big(s,z^{(n-1)}\big)=\sum_{\pi\in\Pi [n-1] } \frac{\partial^{|\pi|} r}{\partial x^{|\pi|}}\big(s, z^{(n-1)}\big) \prod_{B\in\pi}z^{(n-1)}_{|B|}, 
\end{equation} 
\item 
  $\sigma^{(n-1)}:[0,T]\times\real^n\to\real^n$
  is given by 
$$ 
\nonumber 
\sigma^{(n-1)}\big(t,z^{(n-1)}\big)
=\big(\sigma_0^{(n-1)}, \dots, \sigma_{n-1}^{(n-1)}\big)^\intercal
=\big(\sigma^{(n-2)}\big(t,z^{(n-2)}\big),\  \sigma_{n-1}^{(n-1)}\big(t, z^{(n-1)}\big)\big)^\intercal,
$$
where
\begin{equation}
  \label{eq:faa_sigma}
 \sigma^{(n-1)}_{n-1}\big(s,z^{(n-1)}\big)=\sum_{\pi\in\Pi [n-1] } \frac{\partial^{|\pi|} \sigma}{\partial x^{|\pi|}}\big(s, z^{(n-1)}\big) \prod_{B\in\pi}z^{(n-1)}_{|B|}, 
\end{equation}
\item
  $g^{(n-1)}:[0,T]\times\real^n\times\real\to\real^n$
  is given by 
   $$ 
\nonumber 
g^{(n-1)}(t,z^{(n-1)},y)
=\big(g_0^{(n-1)}, \dots, g_{n-1}^{(n-1)}\big)^\intercal
=\big(g^{(n-2)}(t,z^{(n-2)},y),\  g_{n-1}^{(n-1)}(t, z^{(n-1)},y)\big)^\intercal,
$$ 
where
\begin{equation} 
\label{eq:faa_g}
 g^{(n-1)}_{n-1}\big(s,z^{(n-1)},y\big)=\sum_{\pi\in\Pi [n-1] } \frac{\partial^{|\pi|} g}{\partial x^{|\pi|}}\big(s, z^{(n-1)},y\big) \prod_{B\in\pi}z^{(n-1)}_{|B|},
\end{equation}
\end{enumerate} 
where $\Pi [n-1]$ stands for the set of partition of $\{1, \dots, n-1\}$. 
We also assume that $X^{(k)}_t (x)$ is solution of the SDE \eqref{eq:faa_Xn} for
$k = 0,\ldots , n-1$.
\\
 
\noindent
Observe first that \eqref{eq:Zn-1} holds
for $Z^{(0)}=X(x)$ since in this case \eqref{eq:Zn-1} reduces to
\eqref{eq:edsX} with 
$$
r^{(0)}\big(t,z^{(0)}\big)=r\big(t,z^{(0)}\big), 
\quad
\sigma^{(0)}\big(t,z^{(0)}\big)=\sigma\big(t,z^{(0)}\big), \quad
g^{(0)}(t,z^{(0)},y)=g(t,z^{(0)},y).  
$$
Next, regarding $Z^{(1)}$, by \eqref{eq:edsX} and \eqref{eq:Z1} we have 
\begin{align*}
Z^{(1)}(t)  = & 
\left(\begin{array}{l} X_t\\ X^{(1)}_t
\end {array}\right)
=\left(\begin{array}{l} 1\\0\end{array}\right)
+\int_0^t\left(\begin{array}{c} r(s, X_s)\\\frac{\partial r}{\partial x}(s,X_s)X^{(1)}_t\end{array}\right)ds
+\int_0^t\left(\begin{array}{c} \sigma(s, X_s)\\\frac{\partial \sigma}{\partial x}(s,X_s)X^{(1)}_s\end{array}\right)dW_s
\\&
+\int_0^t\int_{-\infty}^\infty \left(\begin{array}{c} g(s, X_{s^-},y) \\ g(s, X_{s^-},y)X^{(1)}_{s^-}\end{array}\right)\big(N(ds,dy)-\nu_s(dy)ds\big)
\end{align*}
which is \eqref{eq:Zn-1} for $Z^{(1)}$ with 
\begin{eqnarray*}
r^{(1)}\big(t,z^{(1)}\big)&=&\left(\begin{array}{c} r(t, z^{(1)}_0)\\ \frac{\partial r}{\partial x}(t,z^{(1)}_0)z^{(1)}_1\end{array}\right)
, \quad
\sigma^{(1)}\big(t,z^{(1)}\big)=\left(\begin{array}{c} \sigma(t, z^{(1)}_0)\\ \frac{\partial \sigma}{\partial x}(t,z^{(1)}_0)z^{(1)}_1\end{array}\right),
\\
g^{(1)}(t,z^{(1)},y)&=&\left(\begin{array}{c} g(t, z^{(1)}_0,y)\\ \frac{\partial g}{\partial x}(t,z^{(1)}_0,y)z^{(1)}_1
\end{array}\right),
\end{eqnarray*} 
corresponding indeed to \eqref{eq:rn}--\eqref{eq:faa_g} in this case. 
We now show that $Z^{(n)}(x)$ solves an SDE
similar to \eqref{eq:Zn-1},
and that $X^{(n)}_t $ is solution to \eqref{eq:faa_Xn} for the index~$n$. 
Since $Z^{(n-1)}(x)$ is solution to \eqref{eq:Zn-1}, by
Theorem~6-29 in \cite{bichteler},
 $\nabla_{z^{(n-1)}}Z^{(n-1)}$ is solution of the $M_{n,n}(\real)$-valued matrix equation
\begin{eqnarray}
\nonumber
\lefteqn{
\nabla_{z^{(n-1)}}Z^{(n-1)}_t(x)=I_{n,n}
+\int_0^t \nabla_{z^{(n-1)}} r^{(n-1)}(s,Z^{(n-1)}_s(x)) \nabla_{z^{(n-1)}}Z^{(n-1)}_s(x) \ ds
}
  \\
\nonumber
&&+\int_0^t \nabla_{z^{(n-1)}}\sigma^{(n-1)}(s,Z^{(n-1)}_s(x))\nabla_{z^{(n-1)}}Z^{(n-1)}_s(x)\ dW_s
\\
\nonumber
&&+\int_0^t\int_{-\infty}^\infty \nabla_{z^{(n-1)}}g^{(n-1)}(s,Z^{(n-1)}_{s^-}(x),y)\nabla_{z^{(n-1)}}Z^{(n-1)}_{s^-}(x)\ \big(N(ds,dy)-\nu_s(dy)ds \big).
\\
\label{eq:Zn-1'}
\end{eqnarray}
With the notation
$\nabla_{z^{(n-1)}}=\Big(\frac{\partial}{\partial x},\frac{\partial}{\partial z_1}, \dots, \frac{\partial}{\partial z_{n-1}}\Big)$, extracting the first column for the matrix equality in \eqref{eq:Zn-1'} we have 
\begin{align} 
\nonumber
& 
  \frac{\partial Z^{(n-1)}_t}{\partial x}(x)=(1,0, \dots, 0)^\intercal
\\
\nonumber
&+\int_0^t \nabla_{z^{(n-1)}} r^{(n-1)}(s,Z^{(n-1)}_s(x)) \frac{\partial Z^{(n-1)}_s}{\partial x} (x)\ ds
+\int_0^t \nabla_{z^{(n-1)}}\sigma^{(n-1)}(s,Z^{(n-1)}_s(x))\frac{\partial Z^{(n-1)}_s}{\partial x}(x)\ dW_s
\\
\label{eq:Zn-1'2}
&+\int_0^t\int_{-\infty}^\infty \nabla_{z^{(n-1)}}g^{(n-1)}(s,Z^{(n-1)}_{s^-}(x),y)\frac{\partial Z^{(n-1)}_{s^-}}{\partial x}(x)\ \big(N(ds,dy)-\nu_s(dy)ds \big).
\end{align} 
Next, for the leftmost 
entry in \eqref{eq:Zn-1'2} we have 
\begin{align} 
\nonumber
 & 
\frac{\partial^n X_t}{\partial x^n}(x)
\\
\nonumber
&= \int_0^t 
\bigg(
\frac{\partial r_{n-1}^{(n-1)}}{\partial x} (s,Z^{(n-1)}_s(x))\frac{\partial X_s}{\partial x} (x)
+\sum_{i=1}^{n-1} \frac{\partial r_{n-1}^{(n-1)}}{\partial z^{(n-1)}_i} (s,Z^{(n-1)}_s(x))\frac{\partial^{i+1} X_s}{\partial x^{i+1}} (x)\bigg)\ ds
\\
\nonumber 
&+\int_0^t \bigg(
\frac{\partial \sigma_{n-1}^{(n-1)}}{\partial x} (s,Z^{(n-1)}_s(x))\frac{\partial X_s}{\partial x} (x)
+\sum_{i=1}^{n-1} \frac{\partial \sigma_{n-1}^{(n-1)}}{\partial z^{(n-1)}_i} (s,Z^{(n-1)}_s(x))\frac{\partial^{i+1} X_s}{\partial x^{i+1}} (x)\bigg)\ \ dW_s
\\
\nonumber
&+\int_0^t\int_{-\infty}^\infty 
\bigg(
\frac{\partial g_{n-1}^{(n-1)}}{\partial x} (s,Z^{(n-1)}_{s^-}(x),y)\frac{\partial X_s}{\partial x} (x)
+\sum_{i=1}^{n-1} \frac{\partial g_{n-1}^{(n-1)}}{\partial z^{(n-1)}_i} (s,Z^{(n-1)}_{s^-}(x),y)\frac{\partial^{i+1}X_s}{\partial x^{i+1}} (x)\bigg)\\
\label{eq:Zn-1''}
& \hspace{10cm}  \big(N(ds,dy)-\nu_s(dy)ds \big).
\end{align} 
Putting together \eqref{eq:Zn-1} and \eqref{eq:Zn-1''}
yields an equation for $Z^{(n)}(x)$ similar to \eqref{eq:Zn-1}, and \eqref{eq:Zn-1''} proves \eqref{eq:faa_Xn} for $X^{(n)}(x)=\frac{\partial^n}{\partial x^n}X(x)$.  
Indeed, from \eqref{eq:Zn-1''} we recover the
expressions of $r_n^{(n)}$, $\sigma_n^{(n)}$ and $g_n^{(n)}$ as in \eqref{eq:faa_r}--\eqref{eq:faa_g}, which achieves the induction. 
For instance, using the expression \eqref{eq:faa_r} of $r_{n-1}^{(n-1)}$, we have 
\begin{eqnarray} 
\nonumber
r_{n}^{(n)}\big(s,z^{(n)}\big) & = & 
z_1^{(n)} \frac{\partial r^{(n-1)}_{n-1}}{\partial x}\big(s,z^{(n-1)}\big)
+ z_{i+1}^{(n)} \sum_{i=1}^{n-1} \frac{\partial r^{(n-1)}_{n-1}}{\partial z_i^{(n-1)}}\big(s,z^{(n-1)}\big) 
  \\
\nonumber
&= & z_1^{(n)} \frac{\partial}{\partial x}\left(\sum_{\pi\in\Pi [n-1]} \frac{\partial^{|\pi|} r}{\partial x^{|\pi|}}\big(s, z^{(n)}\big) \prod_{B\in\pi}z^{(n-1)}_{|B|}\right) 
  \\
\nonumber
& & +\sum_{i=1}^{n-1} z_{i+1}^{(n)}
\frac{\partial}{\partial z_i^{(n-1)}}\left(
\sum_{\pi\in\Pi [n-1]} \frac{\partial^{|\pi|} r}{\partial x^{|\pi|}}\big(s, z^{(n)}\big) \prod_{B\in\pi}z^{(n)}_{|B|}
\right) 
\\
\nonumber
&= & z_1^{(n)} \sum_{\pi\in\Pi [n-1] } \frac{\partial^{|\pi|+1} r}{\partial x^{|\pi|+1}}\big(s, z^{(n)}\big) \left(
\prod_{B\in\pi}z^{(n-1)}_{|B|} \right)
  \\
\nonumber
& & +\sum_{i=1}^{n-1} z_{i+1}^{(n)}
\sum_{\pi\in\Pi [n-1]} \frac{\partial^{|\pi|} r}{\partial x^{|\pi|}}\big(s, z^{(n)}\big) \frac{\partial}{\partial z_i^{(n-1)}}
\prod_{B\in\pi}z^{(n)}_{|B|}.
\\
\label{eq:faa_tech1}
\end{eqnarray} 
 In the second term above we have
\begin{eqnarray}
\nonumber
  \frac{\partial}{\partial z_i^{(n-1)}}
  \prod_{B\in\pi}z^{(n)}_{|B|}
 & = & \frac{\partial}{\partial z_i^{(n-1)}}
 \prod_{j=1}^{n-1} (z^{(n)}_j)^{\#\{B\in \pi : |B|=j\}}
\\
\nonumber
&=&\#\{B\in \pi : |B|=i\}
(z^{(n)}_i)^{\#\{B\in \pi : |B|=i\}-1}
\prod_{\begin{subarray}{c} 1\leq j\leq n\\ j\not=i\end{subarray}} (z^{(n)}_j)^{\#\{B\in \pi : |B|=j\}}.
\\
\label{eq:faa_tech2}
\end{eqnarray}
We note that $\Pi [n]$ consists of partitions $\pi\in\Pi [n-1]$ with
either the addition of $\{n\}$ as a new block,  
or the completion of an existing block by $\{n\}$.
In the latter case, if $\{n\}$ is added to a block of size $j$ the new partition of $\Pi [n]$ obtained
in this way will
have one block of size $j$ less and
one block of size $j+1$ more,
and there are $\#\{B\in \pi : |B|=j\}$ such blocks. 
We conclude that the sums in \eqref{eq:faa_tech1} and \eqref{eq:faa_tech2}
are effectively over $\Pi [n]$,
which yields \eqref{eq:faa_r} for $r_n^{(n)}$, are as follows:
$$
r_{n}^{(n)}\big(s,z^{(n)}\big)=
\sum_{\pi\in\Pi [n]} \frac{\partial^{|\pi|} r}{\partial x^{|\pi|}}\big(s, z^{(n-1)}\big) \prod_{B\in\pi}z^{(n-1)}_{|B|}.
$$
Similar computations yield also \eqref{eq:faa_sigma} and \eqref{eq:faa_g} and achieves the induction proving Proposition~\ref{prop:faa}.  
\end{Proof}

\section{Regularity of stochastic flows} 
\label{sec:regularity}
In this section we consider the solution $(X_t(x))_{t\in [0,T]}$ of SDE \eqref{eq:edsX},
for which Proposition~\ref{prop:faa} gives condition for the
differentiability of the flow $x\mapsto X_t(x)$
up to any order $n$.
 The next Theorem~\ref{theo:moment_derivees}
deals with the integrability of order $q\geq 2$
of the flow derivatives, based on
Corollaries~\ref{prop:Bass_Cranston3}  
and \ref{corol:moment_affineU}, see also Theorem~3.3 in \cite{kunitalevy} 
  which only covers first order differentiability. 
We let $\|f(\cdot )\|_\infty$ denote 
the supremum of functions $f$ on $\rit$.
\begin{theorem} 
\label{theo:moment_derivees}
Let $n\geq 1$ and $q\geq 2$, and assume that \eqref{eq:compensator} and 
(\hyperlink{BGJhyp}{$A_n$}) hold.
Then, for all $k=1,\ldots ,n$ we have 
\begin{equation}
\label{eq:moment_Xk}
\sup_{x\in \real } \ee\bigg[
  \sup_{t\in [0,T]}
  \big| X_t^{(k)}(x)\big|^q\bigg] <\infty.
\end{equation}
\end{theorem}
\begin{Proof}
  Since $X^{(n)}$
  in \eqref{eq:faa_Xn} is
  expressed in terms of $X^{(k)}$ for $k<n$,
  deriving a moment of order $q$ for $X^{(n)}$
  requires to show the existence of
  moments of order $q$ of certain products of the $X^{(k)}$.
  Accordingly, from the H\"older inequality,
  higher moments of every $X^{(k)}$, $k<n$, are required in our argument,
  see \eqref{eq:p*rec3} below. 
 By Proposition~\ref{prop:faa}, 
 $X_t^{(k)}(x):=
 \displaystyle \frac{\partial^k X_t}{\partial x^k}(x)$
solves the SDE
\begin{eqnarray}
\nonumber
\lefteqn{
  dX^{(k)}_t (x) =
  \Big(\sum_{\pi\in\Pi [k]} \frac{\partial^{|\pi|} r}{\partial x^{|\pi|}}\big(x, X_t(x)\big) \prod_{B\in\pi}X^{(|B|)}_t(x)\Big)\ dt
}
\\
\nonumber
&&+\Big( \sum_{\pi\in\Pi [k]} \frac{\partial^{|\pi|} \sigma}{\partial x^{|\pi|}}\big(x, X_t(x)\big) \prod_{B\in\pi}X^{(|B|)}_t(x)\Big)\ dW_t\\
\label{eq:EDS_k}
&&+\int_{-\infty}^\infty\Big(\sum_{\pi\in\Pi [k]} \frac{\partial^{|\pi|}g}{\partial x^{|\pi|}}\big(t, X_t(x),z\big) \prod_{B\in\pi}X^{(|B|)}_t(x)\Big)\ \big(N(dt,dz)-\nu_t(dz)dt\big), \qquad 
\end{eqnarray}
with $X^{(k)}_0 (x)=0$ for $k\geq 2$, where the sum over
$\pi$ runs in the set $\Pi [k]$ of all partitions of $\{1, 2, \dots,k\}$.
In order to prove \eqref{eq:moment_Xk} for $k = 1, \ldots , n$
we shall prove 
\begin{equation}
\label{eq:moment_Xkrec}
\sup_{x \in \real } 
\ee\bigg[
  \sup_{t\in [0,T]}
  \big|X_t^{(k)}(x)\big|^{p_k}\bigg] <\infty \quad 
\end{equation}
by induction, for the order $p_k$ defined by 
\begin{equation}
\label{eq:p*rec}
p_k=q \frac{n!}{k!}. 
\end{equation}
\noindent
For $k=1$, 
 \eqref{eq:EDS_k} reduces to the affine equation 
\begin{equation} 
\nonumber 
dX^{(1)}_t(x) = 
a^{(1)}_t X^{(1)}_t(x)\ dt
+b^{(1)}_t\ dW_t +X^{(1)}_{t^-}(x)\int_{-\infty}^\infty
c^{(1)}_{t^-} (z)\ \big(N(dt,dz)-\nu_t(dz)dt\big),
\end{equation} 
of the form \eqref{eq:affine}, 
with $X^{(1)}_0(x)=1$ and 
$$
a^{(1)}_t=\frac{\partial r}{\partial x}\big(t, X_t(x)\big),
\quad 
b^{(1)}_t=\frac{\partial \sigma}{\partial x}\big(t,X_t(x)\big),
\quad 
c^{(1)}_t(y)= \frac{\partial g}{\partial x}\big(t,X_t (x),y\big), 
$$
 see also Theorem 6-29 in \cite{bichteler}.
 Since 
$(a^{(1)}_t)_{t\in [0,T]} , (b^{(1)}_t)_{t\in [0,T]} \in L^{p_1} ( [0,T], L^\infty (\Omega ))$
under (\hyperlink{BGJhyp}{$A_n$})
 and since Conditions~\eqref{eq:cc_hyp}-\eqref{eq:cc_hyp2}
are satisfied by $c^{(1)}$ with $p=p_1 = n!q$ under \eqref{5.1}, 
 Corollary~\ref{corol:moment_affineU}
 shows that $X^{(1)}_t(x)$ 
 \textcolor{black}{admits a moment of order $p_1$},
 uniformly bounded in $t\in[0,T]$, that is \eqref{eq:moment_Xkrec} holds true for $k=1$. 
 \\
 
  \noindent
Further, we assume that \eqref{eq:moment_Xkrec} holds true with order $p_l$ for
$l = 1, \ldots , k-1$ and we show that it remains true for the rank $l=k$. 
We note that
\eqref{eq:EDS_k} is thus an affine SDE of the form \eqref{eq:affine}
in $X^{(k)}_t(x)$, with the random coefficients
\begin{align*}
& a^{(k)}_{x,t}:= \frac{\partial r}{\partial x}\big(t, X_t(x)\big),
\\
& b^{(k)}_{x,t}:=\frac{\partial \sigma}{\partial x}\big(t, X_t(x)\big),
\\
& c^{(k)}_{x,t}(z):= \frac{\partial g}{\partial x}\big(t, X_t(x),z\big)
\\
& u^{(k)}_{x,t}:=\sum_{\pi\in\Pi [k] \setminus\{1, \dots,k\}} \frac{\partial^{|\pi|} r}{\partial x^{|\pi|}}\big(t, X_t(x)\big) \prod_{B\in\pi}X^{(|B|)}_t(x), 
\\
& v^{(k)}_{x,t}:=\sum_{\pi\in\Pi [k] \setminus\{1, \dots,k\}} \frac{\partial^{|\pi|} \sigma}{\partial x^{|\pi|}}\big(t,X_t(x)\big) \prod_{B\in\pi}X^{(|B|)}_t(x),
\\
& w^{(k)}_{x,t}(z):= \sum_{\pi\in\Pi [k] \setminus\{1, \dots,k\}} \frac{\partial^{|\pi|} g}{\partial x^{|\pi|}}\big(t, X_t(x),z\big) \prod_{B\in\pi}X^{(|B|)}_t(x), 
\end{align*}
as in \eqref{vjdsgf}.
 In order to show that $X^{(k)}_t(x)$ satisfies \eqref{eq:moment_Xkrec}, we shall 
apply Corollary~\ref{corol:moment_affineU}
for every $p_k$ as in \eqref{eq:moment_Xkrec}
to the affine SDE \eqref{eq:EDS_k} written as
\eqref{eq:affine-alpha} and parameterized by the initial condition $x\in\real$, after checking that
$(a^{(k)}_{x,t})_{t\in [0,T]} , (b^{(k)}_{x,t})_{t\in [0,T]} \in L^{p_k} \big( [0,T], L^\infty (\Omega\times\real )\big)$, 
$(u^{(k)}_{x,t})_{t\in [0,T]} , (v^{(k)}_{x,t})_{t\in [0,T]} \in L^{p_k} \big( [0,T] \times \Omega, L^\infty(\real)\big)$,
and the conditions~\eqref{eq:cc_hypU}-\eqref{eq:cc_hypU2} 
for $(c^{(k)}_{x,t}(z))_{t\in [0,T]}$, $(w^{(k)}_{x,t}(z))_{t\in [0,T]}$
 hold, as follows. 
\begin{enumerate}[i)]
\item
The conditions 
$$
\sup_{x\in \real}
\int_0^T\|a_{x,t}^{(k)}\|_\infty^{p_k}\ dt<\infty, \quad
\sup_{x\in \real}
\int_0^T\|b_{x,t}^{(k)}\|_\infty^{p_k}\ dt<\infty,
$$
follow immediately from Assumption~(\hyperlink{BGJhyp}{$A_n$}) on
$\partial r/\partial x$ and $\partial \sigma/\partial x$. 
On the other hand, regarding $c^{(k)}$, the bounds 
$$
\int_0^T\Big\|\int_{-\infty}^\infty |c_{x,t}^{(k)}(z)|^2  \ \nu_t(dz)\Big\|_\infty^{p_k/2}\ dt<\infty
\
\mbox{ and }
\ 
\int_0^T\Big\|\int_{-\infty}^\infty |c_{x,t}^{(k)}(z)|^{p_k}  \ \nu_t(dz)\Big\|_\infty \ dt<\infty
$$
follow from \eqref{5.1} 
since
$\displaystyle
\sup_{x\in \real}
\|c_{x,t}^{(k)}(y)\|_\infty\leq
\displaystyle
\left\|\frac{\partial g}{\partial x}(t,\cdot,y)\right\|_\infty$. 
\item 
    \noindent
 We show that $(u^{(k)}_{x,t})_{t\in [0,T]}, (v^{(k)}_{x,t})_{t\in [0,T]} \in L^{p_k} ( [0,T] \times \Omega )$, uniformly in $x\in \real$. 
Since the coefficients of \eqref{eq:EDS_k} involve finite sums, we can deal with each summand separately using the convexity of $y\in\real\mapsto|y|^{p_k}$. 
For all $\pi\in\Pi [k] \setminus\{1, \dots,k\}$, using
(\hyperlink{BGJhyp}{$A_n$}) and H\"older's inequality yields 
\begin{equation} 
  \label{eq:p*rec3}
  \ee \bigg[ \int_0^T
\Big|\frac{\partial^{|\pi|} r}{\partial x^{|\pi|}}\big(x,X_t(x)\big)\Big|^{p_k} \prod_{B\in\pi}\big|X^{(|B|)}_t(x)\big|^{p_k}\ dt \bigg] 
  \leq C \int_0^T \prod_{B\in\pi}
  \Big(
  \ee\big[\big|X^{(|B|)}_t(x)\big|^{p_k|\pi|}\big]
  \Big)^{1/|\pi|}\ dt.
\end{equation} 
Since $\pi\in\Pi [k]\setminus\{1, \dots,k\}$ has at least two blocks, for $B\in \pi$, we have $|B|<k$ and the induction hypothesis \eqref{eq:moment_Xkrec} applies for $X^{(|B|)}$. 
Additionally, $|\pi|\leq k$, and so 
\begin{equation}
 \label{eq:p*rec4}
p_k|\pi| \leq k p_k
= p_{k-1}.
\end{equation}
As a consequence, by the induction hypothesis \eqref{eq:moment_Xk} applied to each $X^{(|B|)}_t(x)$, we have 
\begin{equation}
 \label{eq:moment_XB}
\sup_{x\in \rit}
\ee\bigg[
  \sup_{t\in[0,T]}
  \big|X^{(|B|)}_t(x)\big|^{p_k|\pi|}\bigg]<\infty
\end{equation}
for each $B\in \pi$ and \eqref{eq:p*rec3} ensures that
$$
\sup_{x\in \real}
\ee \bigg[ \int_0^T |u_{x,t}^{(k)}|^{p_k}\ dt \bigg] <\infty,
$$
 and similarly for $v^{(k)}_{x,t}$ we find:
$$
\sup_{x\in \real}
\ee \left[ \int_0^T |v_{x,t}^{(k)}|^{p_k}\ dt \right] <\infty.
$$
\item Verification of \eqref{eq:cc_hyp} for $w^{(k)}_{x,t}(z)$. 
  Again, since the set $\Pi [k]$ of partitions of $\{1, 2, \dots,k\}$
  is finite and $y\in\rit\mapsto |y|^2$, $y\in\rit\mapsto |y|^{p_k/2}$ are
  both convex functions, we can deal with each summand separately. 
For all $\pi\in\Pi [k] \setminus\{1, \dots,k\}$ we have 
\begin{eqnarray}
\nonumber
\lefteqn{
  \ee \bigg[
    \int_0^T \Big(\int_{-\infty}^\infty \Big|\frac{\partial^{|\pi|} g}{\partial x^{|\pi|}}\big(t, X_t(x),z\big)\Big|^2 \prod_{B\in\pi}\big|X^{(|B|)}_t(x)\big|^2 \ \nu_t(dz)\Big)^{p_k/2} \ dt \bigg] 
}
\\
\nonumber
&=&\ee\bigg[
  \int_0^T\Big(\int_{-\infty}^\infty \Big|\frac{\partial^{|\pi|} g}{\partial x^{|\pi|}}\big(t, X_t(x),z\big)\Big|^2 \nu_t(dz)\Big)^{p_k/2}\prod_{B\in\pi}\big|X^{(|B|)}_t(x)\big|^{p_k} \ dt \bigg]
\\
\nonumber
&\leq&\int_0^T\Big(\int_{-\infty}^\infty \Big\|\frac{\partial^{|\pi|} g}{\partial x^{|\pi|}}\big(t, \cdot,z\big)\Big\|_\infty^2 \nu_t(dz)\Big)^{p_k/2} 
\prod_{B\in\pi}\ee\Big[\big|X^{(|B|)}_t(x)\big|^{p_k|\pi|} \Big]^{1/|\pi|} \ dt\\
\nonumber
&\leq&C\int_0^T\Big(\int_{-\infty}^\infty \Big\|\frac{\partial^{|\pi|} g}{\partial x^{|\pi|}}\big(t, \cdot,z\big)\Big\|_\infty^2 \nu_t(dz)\Big)^{p_k/2} \ dt, 
\end{eqnarray}
where the last bound uses \eqref{eq:moment_XB} for the $X^{(|B|)}_t(x)$ with $|B|<k$
due to the induction hypothesis. 
This final bound is finite under \eqref{5.1}, 
 which ensures that 
$$
\sup_{x\in \real}
\ee\left[
  \int_0^T\Big(\int_{-\infty}^\infty\big|w^{(k)}_{x,t}(z)\big|^2 \ \nu_t(dz)\Big)^{p_k/2} \ dt \right] <\infty. 
$$
 Similarly for \eqref{eq:cc_hyp2}, since $y\in\rit\mapsto |y|^2$ and $y\in\rit\mapsto |y|$ are both convex functions we have 
\begin{eqnarray}
\nonumber
\lefteqn{
  \ee \bigg[
    \int_0^T \int_{-\infty}^\infty \Big|\frac{\partial^{|\pi|} g}{\partial x^{|\pi|}}\big(t, X_t(x),z\big)\Big|^{p_k} \prod_{B\in\pi}\big|X^{(|B|)}_t(x)\big|^{p_k}\ \nu_t(dz) \ dt \bigg] 
}
\\
\nonumber
&=&\ee\bigg[
  \int_0^T \int_{-\infty}^\infty \Big|\frac{\partial^{|\pi|} g}{\partial x^{|\pi|}}\big(t, X_t(x),z\big)\Big|^{p_k} \nu_t(dz) \prod_{B\in\pi}\big|X^{(|B|)}_t(x)\big|^{p_k} \ dt \bigg]
\\
\nonumber
&\leq&\int_0^T \int_{-\infty}^\infty \Big\|\frac{\partial^{|\pi|} g}{\partial x^{|\pi|}}\big(t, \cdot,z\big)\Big\|_\infty^{p_k} \nu_t(dz)  
\prod_{B\in\pi}\ee\Big[\big|X^{(|B|)}_t(x)\big|^{p_k|\pi|} \Big]^{1/|\pi|} \ dt\\
\nonumber
&\leq&C\int_0^T \int_{-\infty}^\infty \Big\|\frac{\partial^{|\pi|} g}{\partial x^{|\pi|}}\big(t, \cdot,z\big)\Big\|_\infty^{p_k} \nu_t(dz) \ dt, 
\end{eqnarray}
using \eqref{eq:moment_XB} for $X^{(|B|)}_t(x)$ with $|B|<n$,
and we conclude to 
$$
\sup_{x\in \real}
\ee\bigg[
  \int_0^T \int_{-\infty}^\infty|w^{(k)}_{x,t}(z)|^{p_k}\ \nu_t(dz) \ dt \bigg] <\infty. 
$$
\end{enumerate}
As a consequence,
Theorem~\ref{theo:moment_affine} can be applied to \eqref{eq:EDS_k}, which yields 
$$
\sup_{t\in[0,T]} \ee\big[\big|X_t^{(k)}(x)\big|^{p_k}\big]<\infty, 
$$
proving the induction hypothesis \eqref{eq:moment_Xkrec} for index $k$ and with order $p_k$ in \eqref{eq:p*rec}. 
In particular, this proves Theorem~\ref{theo:moment_derivees}.
\end{Proof}

\medskip

\noindent
{\em Acknowledgement.} We thank Yufei Zhang for corrections to
Lemma~\ref{prop:Bass_Cranston} and Corollary~\ref{prop:Bass_Cranston2}. 

\footnotesize 

\def\cprime{$'$} \def\polhk#1{\setbox0=\hbox{#1}{\ooalign{\hidewidth
  \lower1.5ex\hbox{`}\hidewidth\crcr\unhbox0}}}
  \def\polhk#1{\setbox0=\hbox{#1}{\ooalign{\hidewidth
  \lower1.5ex\hbox{`}\hidewidth\crcr\unhbox0}}} \def\cprime{$'$}

\end{document}